\documentclass{amsart}
\usepackage[cp1251]{inputenc}
\usepackage[english,russian]{babel}
\usepackage{amsmath}
\usepackage{amssymb}
\usepackage{amsfonts}

\newtheorem{quasitheorem}{quasiTheorem}

\begin{document}

\thispagestyle{empty}

\title {Computational experiments with nilpotent Lie algebra}

\author {V.V.Gorbatsevich}

\maketitle

\section {Introduction}

Studying of properties and structure of the arbitrary finite dimensional Lie algebras mainly reduced to the study of nilpotent Lie algebras. Levi-Maltsev theorem for arbitrary finite-dimensional Lie algebras and Maltsev splitting for solvable Lie algebras result in the need to study the structure of nilpotent Lie algebras and their derivations. But nilpotent Lie algebra are very difficult to study. There are some general statements about their structure and their derivations, but we can not give a clear description of them. Only some special classes of nilpotent Lie algebras (2-nilpotent, filiform Lie algebras, nilradicals of parabolic subalgebras of semisimple Lie algebras et al.) are studied.

In this note we present some experimental results on the general matrix  nilpotent Lie algebras derived by calculations on a computer. Perhaps these results in the future will be able to suggest directions for further research of general nilpotent Lie algebras. It is useful to recall that even Gauss before formulation and proof of his results (in number theory) had done, as shown by the study of its archives, a huge number of calculations. Below we are formulated several results of our calculations in a form of quasitheorems, in particular --- some results about characteristically nilpotent Lie algebras.

All Lie algebras considered below are finite-dimensional and defined over a field $K$ of the characteristic 0, but mostly it will be convenient to assume that $K = \bf R$ (in fact in our calculations Lie algebras will be defined over $ \bf Q $). Recall that by filiform is called a nilpotent Lie algebra with maximal possible nilpotency class $n-1$ (where $ n $ is the dimension of a Lie algebra). A Lie algebra is called characteristically nilpotent if its Lie algebra of derivations is nilpotent  (such Lie algebra is nilpotent automatically). Details of information about nilpotent Lie algebras, used below, can be found in the survey \cite {H}. The calculations used freeware GAP program for Windows (http://www.gap-system.org).

\section {Setting objectives and methodology}

Consider the space $\mathcal L_n (K)$  of $ n $-dimensional Lie algebras over a field $ K $ of zero characteristic. It consists of a set of structure constants $ \{ c_ {ij} ^ k \} $ (forming a structure tensor). If $ \{e_i \} $ is a basis of the Lie algebra $ L $, the  operation of commutation $[e_i, e_j] = c_{ij}^ ke_k$ is defined (uniquely in the given basis) by structure constants $c_{ij}^k $. Since $c_{ij}^k $ are skew-symmetric in the lower indices, the set of structure constants $\{ c_{ij} ^ k \}$ can be regarded as a point in $K^ N$, where $ N = n^2 (n-1)/2$. Moreover, by virtue of the Jacobi identities the set $ \mathcal L_n (K) $ is an algebraic subset of $ K ^ N $. We consider also the space $\mathcal N_n (K)$ --- the space of $n$-dimensional nilpotent Lie algebra, it also is an algebraic subset of $ K^ N$. On these two spaces of Lie algebras the group $GL_n (K)$ acts naturally (this action induced by the natural action of this group on $ K ^ n $). The orbits for these actions consist of all isomorphic Lie algebras (more precisely, of their structural tensors). The Lie algebra $L$ is called rigid if its orbit $\mathcal O_L $ is open (usually here are referring to the Zarisski topology).

The set of Lie subalgebras of given dimension in a fixed Lie algebra is a projective algebraic set (lying in the corresponding Grassmann manifold).

An algebraic variety can be decomposed into a union of irreducible components. We are mainly interested in the variety of nilpotent Lie algebra $ \mathcal N_n$. Description of the irreducible components of
$\mathcal N_n (\bf R)$ is known only for the dimensions $ \le  7$ (see  \cite {H}). In particular, it is known that $ \mathcal N_n (\bf R) $ for $n\le 6$ is irreducible and for each such $n$ has only one open orbit (corresponding to the rigid nilpotent Lie algebra). For $ n =  7$ the space of nilpotent Lie algebras is reducible, it has two irreducible components. Each of these two components contains the orbits of one-parameter families of Lie algebras; the union of all orbits of each family is open in the corresponding irreducible component; all Lie algebras of the two families have two generators. In particular, for $n = 7$ in $\mathcal N_n (\bf R)$ there are no rigid Lie. Rigid nilpotent Lie algebras are known only for $ n \le 6 $ (for such $n$ the number of nilpotent Lie algebra is finite). If $n\ge 8$ the existence of rigid Lie algebras in $\mathcal N_n (\bf R)$ is questioned, but now not possible (despite a number of attempts) to find a rigorous proof of their absence.

Let the $X$  be a topological space (for example, an algebraic variety with the Zarisski topology or --- if $X$ is defined over $ \bf R $ --- with the Euclidean topology). Some property of points of this space is called general (or typical), if it holds for all points of some open and dense subspace in $X$. The points for which this property holds are called general relative of the property. Generally speaking, we can view here any open and dense subset $U$ (corresponding general feature is: "to belong to the subset $U$"). If a property holds for all points of an open set, then it will be called the open one (or unconstrained). For example, the property of nilpotent Lie algebra to have two generators is open. Almost all "computer" part of this article will be devoted to nilpotent Lie algebras with two generators. In particular, we will consider some general (typical) properties of such Lie algebras.

For irreducible components of an algebraic variety a property is general if and only it is open. It follows from the  fact that the closure of the corresponding open subset gives all irreducible component, i.e. this subset will be dense. Therefore, when we study the general properties of irreducible components, we can confine ourselves to the study of open properties.

Here are some simple general properties of Lie algebras $\mathcal N_n (\bf R)$.

1. Equality $ \text {codim} _N ([N, N]) = 2 $. The fact is that for nilpotent Lie algebras always $ \text {codim} _N([N, N])\ge 2$. Consequently, codimension 2 is the  minimum possible and therefore the property of commutant to have codimension 2 is open.

2. The property to have two generators. It is known that the minimum number of generators of a nilpotent Lie algebra is equal to $ \text {codim}_N([N, N])$, and therefore the property 2 is equivalent to the property 1.

3. The equality $\dim Z (N) = 1$. This is obviously a general property, since the dimension of the center of a nilpotent Lie algebra is always positive.

Let $N_m(\bf R)$ be the set of all nilpotent upper-triangular real matrices of order $m$. In $N_m (\bf R)$ there are, as we know from the linear algebra, only a finite number of pairwise non-conjugate (relative to the action, induced by the action of the group $GL_m(\bf R)$) one-dimensional Lie subalgebras. They correspond to the different kinds of Jordan nilpotent matrices $\oplus J_k (0)$ --- the direct sums of Jordan blocks corresponding to the eigenvalue 0. In this case, the "general" matrix is a Jordan block $J_m(0)$ of maximal order. In other words, among all nilpotent matrices $A$ in $N_m (\bf R)$ the matrices which are similar to $J_m(0)$ form an open and dense subset. For them $A^m = 0$, but $A^{m-1} \ne 0$. For all other matrices $A^{m-1} = 0 $, which specifies a proper algebraic subvariety of codimension 1.

Explicit constructions of nilpotent Lie algebras present a lot of  difficulties. We can construct a nilpotent Lie algebra as the quotient Lie algebra of free nilpotent Lie algebra. It is useful to know that free nilpotent Lie algebras are finite-dimensional. However, to find the ideals in them (which gives us by factorization arbitrary nilpotent Lie algebras) is rather inconvenient and difficult. Another way to construct nilpotent Lie algebras --- through a sequence of central extensions. This approach has been used in \cite {L} for studying general nilpotent Lie algebras by a computer. However, in \cite {L} there was a methodological feature --- there  considered only the sequences of one-dimensional (i.e., with one-dimensional kernels) central extensions by selecting the corresponding two-dimensional cocycles with coefficients in $\bf R$ on expandable Lie algebras. Received nilpotent Lie algebras are automatically filiform (because general two-dimensional cocycles were taken, and they are always non-zero one's). To construct the arbitrary general (or typical in the  terminology of \cite {L}) nilpotent Lie algebra central extensions must be considered not only with one-dimensional kernels, but also having an arbitrary dimension. Resulting experimental data in \cite {L} (for example, the fact that there were obtained only characteristically nilpotent Lie algebras) can be attributed to the methodology selected by the author (since the property of characteristical nilpotency is general for filiform Lie algebras --- see, for example, \cite {H}).

We are going to build nilpotent Lie algebras using another special approach --- we will consider Lie subalgebras of the Lie algebra $N_m (\bf R) $ of all real nilpotent matrices of order $m$. The Lie algebra $N_m (\bf R) $is nilpotent and therefore any of its Lie subalgebra is nilpotent too. As is known, any nilpotent Lie algebra has a faithful linear nilpotent representation or, in other words, can be embedded as a Lie subalgebra into Lie algebra $N_m(\bf R)$. Thus we see that in limiting ourselves to Lie subalgebras of $N_m(\bf R) $ for all possible values of $m$, we will not miss any finite-dimensional nilpotent Lie algebra. We note some properties of Lie algebras $N_m(\bf R)$:

1. $N_m (\bf R)$ is the nilpotent Lie algebra of nilpotency class $m-1$.

2. The center $Z = Z (N_m(\bf R))$ of the Lie algebra $N_m (\bf R)$ is one-dimensional. It consists of matrices whose nonzero elements are only the upper-left elements $ a_ {1, m} $.

3. The Lie algebra $ N_2 (\bf R) $ is one-dimensional and Abelian, and $ N_3 (\bf R) $ is three-dimensional. Any subalgebra in $ N_3 (\bf) $, generated by two non-commuting elements, is three-dimensional and coincides with $ N_3 (\bf R) $.

\section {Computing experiments}

The author has made computing of "general" Lie subalgebras of $N_m(\bf R)$ (with some values of $m$), generated by two elements. First, we describe our scheme of computational experiments.

We are interested in general Lie subalgebras of $ N_m(\bf R)$ with two generators. Recall that the property of nilpotent Lie algebra to have two generators is general and therefore it is general for each irreducible component of the space of Lie subalgebras . Let $ X, Y $ be two arbitrary matrix of $N_m (\bf R)$, then by $ N (X, Y) $ we denote the Lie subalgebra of $N_m (\bf R)$, generated by these matrices. Note that the set of Lie subalgebras of fixed dimension in a given Lie algebra is an algebraic subset. In order to produce the general Lie subalgebras, we must take as matrices $X,Y$ the general elements of the $ N_m (\bf R)$. But, as stated above, all general elements here are similar to Jordan block $J_m(0) $. We will consider subalgebras up to similarity, and therefore one of the general matrices $ X, Y $, we can assume equals to $J_m(0)$. Let $X = J_m (0)$ and $Y$ is some general matrix in $ N_m (\bf R) $ (which is similar to $J_m(0)$). With the help of a random selection of elements of this matrix we will get "general" matrix $ Y $.

The calculation procedure was as follows:

1. Set the value of the order of the matrices $m$ (in fact it varied from 4 to 12).

2. Set $X = J_m(0)$ (Jordan block of order $m$).

3. Set the matrix $Y$. To do this, we generate random values for its non-zero elements $y_{ij}$ (with $ i <j $) as random integers (with uniform distribution) in the segment $[- 10,10]$ (i.e. we have 21 possibilities for each nonzero matrix element).

4. Next, use the appropriate procedures of the GAP program.

a). First, we calculate the basis of the Lie algebra $N(X,Y)$ and its dimension.

b). Then the dimensions of the upper central series of the Lie algebra $N(X,Y)$ are found.

c). Next, we calculate the algebra $Der$ of derivations of the Lie algebra $N(X,Y)$ and its dimension. We check whether the Lie algebra $Der$ is nilpotent (or, what is the same, whether the Lie algebra $N (X, Y)$ is characteristically nilpotent).

d). It is possible to do other various additional calculations. In particular, we calculated the commutator algebra of $N (X, Y)$ and checked its characteristical nilpotency (in fact it had never taken place).

e). In conclusion, for the randomly selected ideal of codimension 1 we check whether it will be characteristically nilpotent.

The meaning of this fact is that among the constructed (for $m = 6$) 8-dimensional Lie algebra there are no rigid one's (for which this property of codimension 1 ideals, as is known, must be true (\cite {C}) This is confirmed (but alas, it is not proves even as a quasistatement, which discussed below), the assumption is that not only for the dimension 7, but also for the dimension 8 and above rigid nilpotent Lie does not exist.

It turned out that these Lie ideals are never been characteristically nilpotent. But there are examples of nilpotent Lie algebras, which are characteristically nilpotent and any of their ideals of codimension 1 --- too.

All of these computations for a fixed $m$ carried out several times to make sure that answers which we get in this numerical experiments do not change (it was so in all the dozens of calculations carried out by us). The randomization procedure described above resulted in the fact that every time we receive different (i.e. different as subspace in the space of nilpotent matrices) Lie subalgebras of $ N_m(\bf R)$.

Calculations were carried out over the field $\bf Q $. It would be possible to carry out computations and over finite fields, but the results would need to compare with those obtained for the $\bf Q $. We are not going to do it here.

Here are the results of computational experiments (for brevity a Lie algebra $N(X,Y)$ is denoted below as $N$, and its nilpotency class through $l$):

\smallskip

$ N = 4$: $ \dim N = 4$, $ l = 3$

$ \dim C^i (N) = 4, 2, 1 $

$ \dim Der (N) = 7$; nonnilpotent

$ \dim Der [N, N] = 4$ ; nonnilpotent

\medskip

$ N = 5$: $ \dim N = 6$, $ l = 4$

$ \dim C^i(N)= 6, 4, 3, 1 $

$ \dim Der (N) = 10$; nonnilpotent

$ \dim Der [N, N] = 16$; nonnilpotent

\medskip

$ N = 6 $: $ \dim N = 8 $, $ l = 5$

$ \dim C^i(N) = 8, 6, 5, 3, 1$

$ \dim Der (N) = 12$; nilpotent

$ \dim Der [N, N] = 24$; nonnilpotent

\medskip

$ N = 7 $: $ \dim N = 11 $, $ l = 6$

$ \dim C^i (N) = 11, 9, 8, 6, 3, 1$

$ \dim Der (N) = 18$; nilpotent

$ \dim Der [N, N] = 40$; nonnilpotent

\medskip

$N = 8$: $ \dim N = 14 $, $ l = 7$

$ \dim C^i (N) = 14, 12, 11, 9, 6, 3, 1 $

$ \dim Der (N) = 21$; nilpotent

$ \dim Der [N, N] = 53$; nonnilpotent

\medskip

$ N = 9 $: $ \dim N = 18$, $l = 8$

$ \dim C^i (N) = 18, 16, 15, 13, 10, 6, 3, 1$

$ \dim Der (N) = 27$; nilpotent

$ \dim Der [N, N] = 73$; nonnilpotent

\medskip

$ N = 10 $: $ \dim N = 22 $, $l = 9$

$ \dim C^i (N) = 22, 20, 19, 17, 14, 10, 6, 3, 1$

$ \dim Der (N) = 32$; nilpotent

$ \dim Der [N, N] = 86$; nonnilpotent

\smallskip

Note that for all the above cases in our numerical experiments the dimension (in descending order) of the members of the upper and lower central series (and probably members of these series themselves) coincide with each other (for each fixed $m$). 

It also shows that the class of nilpotency of Lie algebra $ N $ equals to the class of nilpotency of enveloping Lie algebra $ N_m(\bf R)$ (which is equals to $m-1$), i.e. they are the maximum possible for the Lie subalgebras of $ N_n (\bf R) $.

 wecan study also general 3-generated Lie subalgebras of $N_m (\bf R)$. However, the volume of calculations (and dimensions of resulting nilpotent Lie algebras) significantly increased. If $m = 7$ the dimension of the general Lie algebra  $N(X, Y, Z)$ is equal to 16, the dimension of the algebra of its derivations equals to  22, and $\dim C^i = 16, 13, 10, 6, 3, 1 $. If $m=6$, the dimension of the Lie algebra $N(X, Y, Z)$ is equal to 12, $\dim C^i = 12, 9, 6, 3, 1 $, the dimension of the algebra of derivations is equal to 16. All these Lie algebras $ N(X, Y, Z)$ always  are characteristically nilpotent.

\section {Analysis of experimental data}

Let us try to grasp the patterns of the dimensions of the sequences of the upper central series of general Lie subalgebra $N (X, Y)$. Comparing the data for $ m = 4, 5, 6, 7, 8, 9, 10$, we see that the difference between the dimensions (which equal to the dimensions of Lie factor-algebras $\dim C^i/C^{i- 1}$) begin in all cases the same: 2, 1, 2, 3 ... and end the same way, too: ... 3, 2, 1.

Let us consider the initial values of the sequence of differences and compare them with the corresponding values for a free nilpotent Lie algebra $ f_2(l)$ of nilpotency class l. They are, as it is well known (see for example \cite {B}, Chapter II, 2.10), as follows: $2, 1, 2, 3, 6 \dots $. We see that for $N (X, Y)$  we have the same sequence at the beginning. Therefore, the initial terms of the sequence $ \dim C^i/C^{i-1} $ are the same as for the free Lie algebra with two generators. But for a free Lie algebra the dimensions of members of upper central series are monotonically increasing, but for $ N (X, Y) $ this is not so --- at the end of the sequence the dimensions of its members decrease. This is similar to the behavior of the members of the central row for $ N_m (\bf R) $, which ends with the values of dimensions 3,1.

For $\dim N (X, Y)$ with $m$ started from 2 we have such segment of values: 1, 3, 4, 6, 8, 11, 14, 18, 22, 26, 35, 42. For values $4 \le  m \le 10 $, there is an exact formula $ \dim N (X, Y) = [m(m + 1)/5]$ (integer part of $m(m + 1)/$ 5). But for other values of $m$ (for $m=3$ or $m>10$) this formula is not correct.

We see that for small $m$ the dimension of the general 2-generated subalgebra of $ N_m (\bf R)$ is approximately 2/5 of the dimension of the whole Lie algebra $N_m(\bf R)$ (which, obviously, equals to $m(m-1)/2$).

\section {Conclusion}

We introduce the concept of quasitheorem  to refer to the results, for which there are fairly detailed explanation by computational experiments.

In the article the following results are "proved" (in the computer sense):

\begin {quasitheorem} For general $X,Y \in N_m(\bf R)$ generated subalgebras $N(X, Y)$ for $m \ge 6$  are characteristically nilpotent.

General ideals of codimension 1 in $N(X, Y)$ are not characteristically nilpotent.

The Lie algebras $N(X,Y)$ are not rigid.

\end {quasitheorem}

The fact that for $m <6$ Lie algebras $ N (X, Y) $ are not characteristically nilpotent, can be explained very simply --- in this case $\dim N (X, Y) \le 6$, but characteristically nilpotent Lie algebra of dimension $\le 6$ does not exist.

Our experiments also studied the properties of characteristic nilpotence for some ideals (as mentioned above).

\begin {quasitheorem}
For general Lie algebra $N(X,Y)$, generated by $X,Y \in N_n (\bf R)$, the Lie algebra $[N (X, Y), N (X, Y)]$ is not characteristically nilpotent.

\end {quasitheorem}

Note that in \cite {L} there is an example (of dimension 18 with very high values of the structure constants; for example, one of them is equals to 2880) of a nilpotent Lie algebra for which the commutator algebra is also characteristically nilpotent. Moreover, it turned out there that for $n \ge 13 $ for general $n$-dimensional nilpotent Lie algebra commutator algebra was always characteristically nilpotent. However, recall that the notion of general Lie algebras there is very specific --- all of them are filiform Lie algebras.

\begin {quasitheorem}
The sequences of dimensions of the members of the upper and lower central series of the Lie algebra $N(X, Y)$ from quasitheorem 1 are identical and are such that their successive difference (equal to $\dim C^i/ C^{i-1})$) are, respectively, $2, 1, 2, 3, \dots, 3, 2, 1 $. The beginning of this sequence --- like for the free (nilpotent) Lie algebra with two generators, as well as in the end --- as for Lie algebra $N_m (\bf R)$.
\end {quasitheorem}

One could also calculate the dimension of members of the derived series for general Lie algebra $N(X,Y)$. For example, for $m =  9$ we get the derived series dimensions 18,16,8.

It would be interesting to find out to what extent the resulting for random matrices $Y$ the the Lie algebras $N(X, Y)$ are isomorphic. The simplest methods of testing isomorphism of Lie algebras over $\bf Q $ and using GAP did not lead to a specific outcome. Dimension series of classical ideals, algebras of derivations of the algebra Lie and it commutant appeared at each experiment, are the same (for a fixed $m$).

Next our quasiresult deals with the question of the existence of new classes of characteristically nilpotent Lie algebras. Initially, they were built with great difficulty, but then it was proved that the property of characteristical nilpotency is general for filiform Lie algebras. Recall that the filiform Lie algebras have two generators. Whether there are characteristically nilpotent Lie algebra with two generators, which are not  filiform? Direct sums of characteristically nilpotent Lie are also characteristically nilpotent, but they are no longer filiform and they have more than two generators. The above procedure of our computations gives us a sequence of nonfiliform nilpotent Lie algebras $N_k = N(X_k,Y_k)$, which are characteristically nilpotent and have two generators, and their dimensions increase indefinitely. This proves

\begin {quasitheorem} There is an infinite sequence of characteristically nilpotent, but not filiform, Lie algebras with two generators, whose dimensions increase indefinitely.
\end {quasitheorem}

There are many reasons to believe that quasitheorems  stated above are correct in the usual mathematical sense, i.e. they can be proved by rigorous deductive method.

Note also that the subalgebra of the form $ N (X, Y) $ in $ N_m (\bf R) $ can be regarded not only as a general 2-generated Lie subalgebra, but also as a maximal 2-generated Lie subalgebras in $N_m (\bf R)$ or maximum available within $N_m $. For these classes of nilpotent Lie some quasitheorems, similar to the above, may be formulated.

\section {Technical details}

The computation time required to prove our quasitheorems, increases rapidly with increasing of the order of matrices $ m $. If $ m = $ 9 counting time for the algebra of differentiations --- about 5-10 minutes, loading the processor by about 12 percent, while for $ m = $ 10 --- about an hour; the calculation itself of the Lie algebra $N (X, Y)$ is nearly instantaneous. It was used a computer with a frequency of 3.7 GHz processor. and 16 MB of RAM (actually, the calculations used in no more than 3 MB. RAM).

Here is the main part of our program for GAP. We set there $m = 10$, but this value can be changed.

m:=10;

J:=NullMat(m,m);;

p:=m-1;

for i in [1..p] do

J[i][i+1]:=1;

od;

Jm:=LieObject(J);

A:=NullMat(m,m);;

for i in [1..p] do

q:=i+1;

for j in [q..m] do

A[i][j]:=Random(-10,10);

od;

od;

Am:=LieObject(A);

L:=LieAlgebra(Rationals, [Am,Jm]);

LieLowerCentralSeries(L);

D:=Derivations(Basis(L));

IsLieNilpotent(D);

\end {document}